\documentclass[11pt]{article}
\usepackage[a4paper]{anysize}\marginsize{3.1cm}{3.1cm}{0.9cm}{3cm}
   \pdfpagewidth=\paperwidth \pdfpageheight=\paperheight
\usepackage[latin1]{inputenc}
\usepackage{graphicx}

\pdfpagewidth=\paperwidth \pdfpageheight=\paperheight
\usepackage{amssymb}
\makeatletter
\renewcommand\@seccntformat[1]{\csname the#1\endcsname.\enspace}
\renewcommand\@begintheorem[2]{\trivlist\item[\hskip\labelsep{\bfseries#1 #2.}]\it}
\renewcommand\@opargbegintheorem[3]{\trivlist\item[\hskip\labelsep{\bfseries#1 #2}] {\bfseries(#3).}\enspace\it\ignorespaces}
\makeatother

\newtheorem{thm}{Theorem}[section]

\newtheorem{corollary}[thm]{Corollary}
\newtheorem{conjecture}[thm]{Conjecture}

\newcommand\mkthm[2]{\newenvironment{#1}{\begin{#2}\rm}{\end{#2}}}
 \mkthm{definition}{tdefinition}
         \mkthm{remark}{tremark}
       \mkthm{example}{texample}
     \mkthm{question}{tquestion}

\renewcommand\tilde{\widetilde}

\newcommand\C{\mathbb C}

\newcommand\R{\mathbb R}
\newcommand\Z{\mathbb Z}
\newcommand\F{\mathbb F}
\newcommand\IP{\mathbb P}
\renewcommand\P{\mathbb P}

\newcommand\be[1][@{\;}r@{\;}c@{\;}l@{\;}l@{\;}]{$$\everymath{\displaystyle}\renewcommand\arraystretch{1.2}\begin{array}{#1}}
\newcommand\ee{\end{array}$$}

\newcommand\compact{\itemsep=0cm \parskip=0cm}

\newcommand\HM{F_{{\mathrm {HM}}}}

\begin{document}

   \title{Wolf Barth (1942--2016)}
   \author{Thomas Bauer, Klaus Hulek, Slawomir Rams, Alessandra Sarti, Tomasz Szemberg}
   \date{}
   \maketitle
   \thispagestyle{empty}



\section{Life}

\subsection{Youth}
Wolf Paul Barth was born on 20th  October 1942 in Wernigerode, a small town  in the  Harz mountains. This is where his family, originally living in Nuremberg,
had sought refuge from the bombing raids. In 1943 the family moved to Simmelsdorf, which is close to Nuremberg, but  less vulnerable to attacks.
In 1945 Barth's brother Hannes was born.
The family, the father was a town employee, then moved back to Nuremberg and Wolf Barth's home was close to the famous Nuremberg castle -- his playgrounds were the ruins of Nuremberg.

In 1961 Wolf Barth successfully completed the Hans Sachs Gymnasium
in Nuremberg.
He was an excellent student (except in sports), and
he then chose to study mathematics and physics at the nearby Universit\"at Erlangen, officially called Friedrich-Alexander Universit\"at Erlangen-N\"urnberg.

\subsection{Studies and early academic career}
As it was the standard at the time, Wolf Barth enroled for the {\em Staatsexamen}, the German high school teacher's examination. At this time, Reinhold Remmert was a professor of mathematics in
Erlangen and Wolf Barth soon felt attracted to this area of mathematics. So, when Remmert moved to G\"ottingen in 1963, Barth went with him.
G\"ottingen had at that time started to recover to some extent after its losses in the Nazi period  and Hans Grauert was one on the leading figures who
attracted many talented mathematicians.
And indeed, it was  Grauert, who became Barth's second academic teacher and a figure  who influenced his mathematical thinking greatly.
Wolf Barth not only completed his studies in G\"ottingen, but in May 1967  he also obtained his PhD with a thesis on   {\em Einige Eigenschaften analytischer Mengen in kompakten komplexen Mannigfaltigkeiten} (Some properties of analytic sets in compact complex manifolds).
In 1967 Remmert moved again, this time to take up a chair at the Westf\"alische Wilhelms-Universit\"at M\"unster where Heinrich Behnke had been a professor, and whose successor Remmert became.
M\"unster, the home of the Behnke school, was at that time one of the international centres of the school of several complex variables.
Again, Barth followed Remmert and worked for the next two years in M\"unster.  These were lively times at German universities (the '68 revolt) and Barth became interested in and attracted to the ideas of the 1968
generation.

Wolf Barth went to spend the academic year 1969/70  as a visiting lecturer  at MIT in Cambridge, USA, a period where he encountered many new faces and ideas.
In 1971, a year after returning from the USA, he obtained his
Habilitation in M\"unster, where he also became a professor.
In 1972 Wolf Barth took up a chair
at Rijksuniversiteit Leiden in the Netherlands. At the time he was the youngest full professor in Leiden. It was here that he started his close collaboration with Antonius Van de Ven -- an encounter which would
greatly influence his future mathematical interests and career, as we will  outline below. Barth  married his wife Regina in 1972.

\subsection{The Erlangen period}
In 1976 Wolf Barth was offered a professorship at Universit\"at Erlangen, his original Alma Mater.
Erlangen is of course well known to mathematicians, not least  through the Klein programme, made famous by the inaugural lecture of Felix Klein, who was a professor there from 1872 to 1875.
Another famous professor had been Max Noether, one of the 19th century masters of algebraic geometry -- the field of mathematics which had become Wolf Barth's own area of research.
Also, his daughter Emmy Noether had been a student in Erlangen where she wrote her thesis in invariant theory under the guidance of Paul Gordan.

This call also reflected the international reputation which Barth had gained in such a short period on the strength of his research; he was offered the very chair
which had been established for Max Noether and which had been held by such prominent mathematicians as  Heinrich Tietze, Johannes Radon,
Wolfgang Krull  and Georg N\"obeling.
This in itself was a great honour. But it was also the chance to move back to Franconia, his home region, which
made this offer irresistible. Wolf Barth truly loved his Franconia and in fact he stayed there for the rest of his life.  In 1977 his son Matthias was born and in 1981 his daughter Ursula followed.

In Erlangen, Barth
took on multiple duties which came with the position of {\em Lehrstuhl\-inhaber} (full professor). It should be emphasized that
teaching was not simply a duty for him, he took it very seriously
and was very conscientious -- we will come back to this aspect of his mathematical life later.
As one of the
full professors Barth  was also constantly involved in running the Institute of Mathematics  (which later became part of the Department of Mathematics). In particular,  he was   Dean of
Naturwissenschaftliche Fakult\"at I from 1981 to 1983.

Barth's arrival in Erlangen created a very lively research atmosphere. In fact, there were  two active research seminars related to algebraic geometry in Erlangen at the time: the seminar run by Wolf Barth and
the seminar organized by Wulf-Dieter Geyer and Herbert Lange.  Often one of the seminars was devoted to lectures on ongoing or recently completed research  with many guests and visitors as speakers,
whereas for the other seminar a subject was chosen to be presented in much detail by the participants of the seminar in turn.
In the early 1980s preparations were made to apply to Deutsche Forschungsgemeinschaft (DFG) for a major research project in complex algebraic geometry.
Several research groups collaborated to submit a proposal for a Priority Programme (Schwerpunktprogramm), including the newly
established algebraic geometry group in Bayreuth, where Michael Schneider had taken up a chair, as well as the group led by Otto Forster in M\"unchen. The application was successful and from 1985 onwards
the DFG Schwerpunktprogramm {\em Komplexe Mannigfaltigkeiten} became an important factor of the German research activities in algebraic geometry. Wolf Barth became the coordinator of the programme and
oversaw two successful applications for an extension of the programme. This was a very active and fruitful period for the researchers involved and it brought many visitors, both short and long term,
to the participating groups in Germany. In particular, the extra funding allowed Wolf Barth to invite students from various countries (Mexico, Israel, Italy, Poland to name a few) to work under his guidance
 -- we will discuss some results  of his efforts to support development of research groups abroad  in the last section of this account.

For many years Barth was a regular organizer of Oberwolfach meetings.
The meeting
{\em Mehrere komplexe Ver\"anderliche},  originally organized  by Grauert, Remmert and Stein, is one of the oldest series of Oberwolfach conferences.
In 1982 Wolf Barth was asked to replace Karl Stein as one of the organizers. In this capacity, together with Grauert and Remmert, he was in charge  of  the biannual meetings until 1994, when the organization of the series was taken over by Demailly, Hulek and Peternell.
Apart from the regular meetings in this series, he often organized Oberwolfach workshops on more specialized topics, typically with Van de Ven as a co-organizer.
In spite of these numerous obligations, Barth also served the mathematical community as an editor.
Most importantly, he and Wolf von Wahl were joint Editors-in-Chief of Mathematische Zeitschrift from 1984 to 1990.

Wolf Barth retired on 1st  April 2011. Quite tellingly, the lecture  he gave during the farewell conference  ``Groups and Algebraic Geometry",
organized by the Emmy-Noether Zentrum to honour his achievements, was called ``99 Semester Mathematik". Even after Barth's retirement,  his former students sometimes obtained
mails with his comments on their recent papers,
and his teaching manuscripts were further  available from the webpage of Universit\"at Erlangen and were widely used by students. Wolf Barth died on December 30, 2016 in Nuremberg.


\section{Wolf Barth -- Research}


\subsection{Barth-Lefschetz theorems}

Wolf Barth's field of research was complex algebraic geometry and his original approach was strongly influenced by the German school of several complex variables, also
known as the Behnke school, named after its founder Heinrich Behnke, and later led by Karl Stein, Hans Grauert and Reinhold Remmert.
Barth became first famous through his work on  the topology of subvarieties of projective space $\mathbb P^N$.
The starting point of this work is the  celebrated Lefschetz hyperplane theorem, which compares the topology of a projective manifold to that of a hyperplane section:
\begin{thm}[Lefschetz hyperplane theorem]\label{teo:NoetherLefschetz}
Let $X \subset {\mathbb P}^N$ be a subvariety of dimension $k$ and let $Y=X \cap H$ be a hyperplane section such that $U=X \setminus Y$ is smooth. Then the restriction map  $H^k(X) \to
H^k(Y)$ in (singular) cohomology is an isomorphism for $k< N-1$ and injective for $k=N-1$.
\end{thm}

Lefschetz proved this result using the by now famous technique of Lefschetz pencils. Another approach was later developed by Andreotti and Fraenkel using Morse theory.
In his paper \cite[p. 952]{B2} Barth generalized the Lefschetz theorem in the following way.
\begin{thm}[Barth]
Let $X,Z \subset {\mathbb P}^N$ be manifolds of dimension $n$ and $m$ respectively, such that $2n\geq N+s$ and $n+m\geq N+r$.
Then $H^k(X) \to H^k(X \cap Z)$ is an isomorphism for $k \leq \min\{r-1,s\}$.
\end{thm}
Barth's proof is analytic. It is based on his earlier work on extending meromorphic functions \cite{B1} and uses $q$-convexity, sheaf cohomology and the Leray spectral sequence.
Barth and Larsen  \cite{BL} further extended the work of Lefschetz from cohomology  to homotopy groups.
Again, their approach uses analytic tools, such as distance functions in projective space and pseudo-concavity of tubular neighborhoods.
In particular, they proved \cite[Theorem I]{BL}:
\begin{thm}[Barth-Larsen]\label{teo:BarthLarsen}
Let $X\subset {\mathbb P}^N$ be smooth of dimension $n$ with $2n\geq N+1$.
Then $X$ is simply connected, i.e. $\pi_1(X)=0$.
\end{thm}
This was later  strengthened further by Larsen \cite{L}. Indeed, the concept of Barth-Lefschetz theorems became a well known term in the literature. The work of Barth and Larsen also
influenced Fulton and Hansen \cite{FH} when they proved their famous connectedness theorem. This circle of ideas was further taken up in the work of Badescu, Sommese and Debarre leading,
among other things, to Barth-Lefschetz theorems in products of projective spaces and, more generally, homogeneous varieties.


\subsection{Vector bundles}\label{subsec:vectorbundles}

After moving to Leiden, Wolf Barth started collaborating with Van de Ven. It was at this time that Hartshorne \cite[p. 1017]{H} formulated his well-known conjecture on complete
intersections, which, in its simple form, can be stated as
\begin{conjecture}[Hartshorne]
Let $X \subset \mathbb P^N$ be a smooth manifold of codimension $2$. If $N \geq 7$ then $X$ is a complete intersection $X=S_1 \cap S_2$ of two hypersurfaces.
\end{conjecture}
This conjecture was motivated, on the one hand, by a number of examples in small dimension and, on the other hand, by the fact that projective  submanifolds in this range are
subject to a number of strong topological constraints.
Indeed, by the Lefschetz theorem, the Picard group, i.e.~the group of line bundles on $X$, is a free abelian group of rank one, more precisely, the restriction
${\mathrm {Pic}}(\mathbb P^N) \to  {\mathrm {Pic}}(X)$ is an isomorphism. Further, by the Barth-Larsen Theorem \ref{teo:BarthLarsen}, $X$ is simply connected.
This conjecture of Hartshorne's can be rephrased  as a statement about vector bundles on projective space.
The connection is via the {\em Serre construction}, which establishes a relationship between locally complete intersections of codimension $2$ and rank $2$ vector bundles,
generalizing the well-known correspondence between divisors and line bundles.
Applied to projective space  $\mathbb P^N$ of dimension $N\geq 3$ this says the following: Assume that $X \subset \mathbb P^N$  is a codimension $2$ manifold which is
{\em sub-canonical}. Then there exists a rank $2$ vector bundle $E$ on $\mathbb P^N$ and a section $s \in H^0(\mathbb P^N,E)$ such that $X=\{s=0\}$.
Here sub-canonical means that the
 canonical line bundle $K_X= \det T_X$, the determinant of the cotangent bundle,
 is the restriction of  a line bundle on $\mathbb P^N$, i.e.~it is of the form $K_X=\mathcal O_X(k)= \mathcal O_{\mathbb P^N}(k)|_X$ where
 $\mathcal O_{\mathbb P^N}(k)=  \mathcal O_{\mathbb P^N}(1)^{\otimes k}$ and  $\mathcal O_{\mathbb P^N}(1)$ is the hyperplane bundle on $\mathbb P^N$.
 Since the canonical bundle can, by the adjunction formula, be written as $K_X= \det N_{X/ \mathbb P^N} \otimes K_{\mathbb P^N}$, where $ N_{X/ \mathbb P^N}$ is the
 normal bundle of $X$ in $\mathbb P^N$, Serre's construction can also be rephrased in more geometric terms as follows: if the determinant of the normal bundle of a codimension $2$ submanifold
 $X$ of $\mathbb P^N$ can be extended to the surrounding projective space, then so can the normal bundle itself.  Now it follows from the fact that the restriction map defines an isomorphism on the Picard groups,
 that every co\-dimen\-sion $2$
 manifold in $\mathbb P^N$ is sub-canonical, provided $N \geq 5$. In particular, the Serre construction can be applied in this case and we obtain that  $X=\{s=0\}$ for some
 section~$s$ of a rank $2$ bundle $E$ on $\mathbb P^N$. Using that the vector bundle $E$ is uniquely determined, one can argue that $X$ is a
 complete intersection if and only if $E$ is a decomposable rank $2$ bundle, i.e.~a sum of two line bundles. Hence the Hartshorne conjecture can be restated in terms of
 vector bundles as
 \begin{conjecture}
If $N \geq 7$, then every rank $2$ bundle E on  $\mathbb P^N$ decomposes, i.e.~is the sum of two line bundles.
\end{conjecture}
It should be noted that there are also no known indecomposable rank $2$ bundles (apart from the Tango bundle on $\mathbb P^5$ in characteristic $2$) for $N=5,6$.
Barth and Van de Ven \cite[Theorem I]{BVdV} proved the asymptotic version of this conjecture. To describe this, we fix a sequence of linear embeddings $i_N: \mathbb P^N \to \mathbb P^{N+1}$. We say that
a vector bundle $E$ on $\mathbb P^N$ extends to $\mathbb P^{N+1}$ if there exists a vector bundle $E'$ on $\mathbb P^{N+1}$ such that $E=i_N^*(E')$. Similarly, we say that
a submanifold $X \subset \mathbb P^N$ extends to $\mathbb P^{N+1}$ (as a submanifold)  if there exists a submanifold $X' \subset  \mathbb P^{N+1}$ with
$i_N(X)=i_N(\mathbb P^N) \cap X'$.
\begin{thm}[Babylonian Tower Theorem, Barth -- Van de Ven]
A rank $2$ vector bundle $E$ on $\mathbb P^N$ which extends to $\mathbb P^{M}$ for all $M \geq N$, splits into the sum of two line bundles.
\end{thm}
By the Serre construction this implies the
\begin{corollary}
A smooth codimension $2$ submanifold $X \subset \mathbb P^N$ which extends, as a submanifold, to $\mathbb P^{M}$, for all $M \geq N$, is a complete intersection of two hypersurfaces.
\end{corollary}
The paper by Barth and Van de Ven \cite{BVdV} contains
many techniques which later became standard tools in the study of vector bundles, such as the idea to investigate vector bundles on $\mathbb P^N$ by studying their
restriction to lines $L \subset \mathbb P^N$. In  their article, Barth and Van de Ven also rediscover the Serre construction, as did Horrocks, Hartshorne and Grauert and M\"ulich
on different occasions.

One  area where Wolf Barth's influence was, and is, extraordinary is the classification of vector bundles. Classifying vector bundles means, as is typically the case with classification problems in algebraic geometry, that one has to construct a  {\em moduli space}, whose geometric properties one aims to understand. To be able to construct moduli spaces it is necessary to restrict  to {\em stable} objects. As suitable notions of stability one can use
 {\em Mumford-Takemoto stability} or, alternatively, {\em Gieseker stability}.  A vector bundle $F$ of rank $r$  on $\mathbb P^N$ is called {\em Mumford-Takemoto semi-stable} if for
every non-trivial coherent subsheaf $E \subset F$ the inequality
$$
c_1(E)/ {\mathrm {rank } (E)}   \leq  c_1(F)/  {\mathrm {rank }}( F)
$$
holds. Here $c_1=c_1(E) \in H^2(\mathbb P^N, \mathbb Z) = \mathbb Z$
denotes the first Chern class  (which one can interpret as an integer) of a coherent sheaf  $E$. The bundle $F$ is called (properly) stable if strict inequality always holds.
Stable vector bundles are subject to topological constraints: Schwarzenberger \cite[Theorem 7]{Sch} showed that the Chern classes $c_i, i=1,2$, of a stable rank $2$ bundle $F$ on $\mathbb P^2$ satisfy the inequality $c_1^2 < 4c_2$.  In \cite[Corollary 2]{B3}  Barth proved that
this inequality also holds  for stable bundles on any $\mathbb P^n, n\geq 2$. In this paper Barth also demonstrated
the power of the, afterward universally used, method of studying a vector bundle on $\mathbb P^N$ by its restriction to lines.
If $F$ is a rank $2$ vector bundle then one can, after twisting with a suitable line bundle, assume that its first Chern class $c_1(F) \in \{0,-1\}$.
If $c_1(F)=0$ then, by the Grauert-M\"ulich theorem \cite[\S6, Satz 2]{GM}, later
extended to higher rank by Spindler, a stable rank $2$ bundle $F$  splits as the sum of two trivial line bundles on a general line $L$, namely
$F|_L= \mathcal O_{L} \oplus \mathcal O_{L} $. The lines where this is not the case, i.e.~where  $F|_L= \mathcal O_{L}(k) \oplus \mathcal O_{L}(-k) $ for some $k>0$
define a divisor, and hence a curve $C=  C(F)$ in the dual projective plane. These lines were called {\em jumping lines}\label{p:jumping lines} by Barth and $C(F)$ was called the
{\em curve of jumping lines}. Studying the curves of jumping lines and related objects became a central theme of the theory of vector bundles in the years after Barth's paper.

Constructing moduli spaces  and understanding their properties is a crucial question in any classification problem in algebraic geometry. The fundamental theory
for moduli spaces of vector bundles on higher-dimensional varieties  was developed by Maruyama \cite{M1}, \cite{M2} and Gieseker \cite{Gi} (for surfaces).
It was Barth in \cite{B4} who gave a beautiful concrete construction
of the moduli spaces $M^2_{0,c_2}(\mathbb P^2)$ of rank $2$ bundles on $\mathbb P^2$ with even first Chern class. The main geometric result is that every stable rank $2$ vector bundle $F$
on  $\mathbb P^2$ with even first Chern class, which we can assume to be $c_1(F)=0$, can be reconstructed from its curve of jumping lines $C(F)$, which is a
 plane curve of degree $n=c_2(F)$, together with an ineffective theta characteristic $\theta$, i.e. a root of the canonical bundle without sections, on $C(F)$.  The other important contribution of Barth in this paper is, that he
 made, for the first time, systematic
 use of monads in the study of moduli of vector bundles. The concept of monad was introduced  by Horrocks~\cite{Ho}, who considered monads as the elementary building blocks for
 constructing vector bundles.  A {\em monad} is simply a complex
 $$
 A \stackrel{f} \to B \stackrel{g}  \to  C
 $$
where $A,B$ and $C$ are vector bundles, $f$ is a monomorphism of vector bundles, $g$ is an epimorphism such that $g \circ f=0$. The cohomology of this complex
$$
F= \mathrm{ker}(g)/ \mathrm{im}(f)
$$
is then a vector bundle. The advantage of monads is that one can often take $A,B$ and $C$ to be very simple bundles and in this way reduce a moduli problem to a linear algebra question.
For stable rank $2$  vector bundles with $c_1(F)=0$ and $c_2(F)=n$, Barth showed that every such bundle can be realized as the cohomology of a monad of the form
$$
n{\mathcal O}_{\mathbb P^2} \stackrel{a}  \longrightarrow n {\Omega}^1_{\mathbb P^2}(1) \stackrel{c}  \longrightarrow (n-2){\mathcal O}_{\mathbb P^2}(1)
$$
where
${\Omega}^1_{\mathbb P^2}$ is the cotangent
bundle on $\mathbb P^2$.
Moreover, for a given sheaf $E$, the sheaf  $E(k)= E \otimes \mathcal O_{\mathbb P^2}(k)$ denotes the twist of $E$ by $\mathcal O_{\mathbb P^2}(k)$,
and $nE$ is the $n$-fold direct sum of $E$.
The crucial ingredient here is the construction of the map $a$ (which in turn determines the map $c$). It is given by the
middle part $H^1(\mathbb P^2,F(-2))\otimes H^0(\mathbb P^2, \mathcal  O_{\mathbb P^2}(1)) \to H^1(\mathbb P^2,F(-1))$ of the cohomology
module $\bigoplus_k H^1(\mathbb P^2, F(k))$
over the homogeneous coordinate ring $\bigoplus_k H^0(\mathbb P^2, \mathcal O_{\mathbb P^2}(k))$.
Via canonical identifications one can then view $a$ as a {\em net of quadrics} in $\mathbb P^2$ and this net of quadrics determines the pair $(C(F),\theta)$ and vice
versa. In particular, $C(F)$ is the discriminant of the net of quadrics and the theta-characteristic encodes how the net can be reconstructed from its discriminant.

Having thus reduced  the classification problem of stable rank $2$ bundles with even first Chern class to a linear algebra problem, Barth was able to show in \cite[p. 83]{B4}:
\begin{thm}[Barth]
The moduli spaces $M^2_{0,n}(\mathbb P^2)$ of stable rank $2$ bundles with first Chern class $c_1(F)=0$ and second Chern class $c_2(F)=n$ are irreducible
rational manifolds of dimension $4n-3$.
\end{thm}
It is worth mentioning that at that time mathematicians still thought it conceivable that all moduli spaces in algebraic geometry are (uni-)rational. It was only at the beginning of the 1980's that Frei\-tag, Tai and Mumford proved that the moduli space $\mathcal A_g$ of principally polarized abelian varieties of dimension $g$ is of general type for $g\geq 7$. This was the first time that such
a phenomenon was observed. As it turned out, Barth's proof gave unirationality rather than rationality of $M^2_{0,n}(\mathbb P^2)$ and the proof of rationality was  finally
completed by Maruyama \cite{M3}.
Barth and Hulek then  studied monads more systematically in their paper \cite{BH}. At the same time, Beilinson developed his very general approach, leading to the
Beilinson spectral sequence \cite{Be}, thereby linking the classification of vector bundles to derived categories.
The geometry of moduli spaces on projective spaces was further pursued and advanced by many authors including Drezet - Le Potier \cite{DLP},
Ellingsrud \cite{Ell}, Forster - Hirschowitz - Schneider \cite{FHS}, Hirschowitz \cite{Hi},  Hulek \cite{Hul}, Le Potier \cite{LP},  Str{\o}mme \cite{Stro},   and many others, cf. also the book by Okonek, Schneider and Spindler  \cite{OSS}.
It is now known that all moduli spaces  $M^r_{c_1,c_2}(\mathbb P^2)$ of stable vector bundles on the projective plane are irreducible and rational. This is, however no longer true on higher-dimensional projective spaces, starting with $\mathbb P^3$ \cite[Section 8]{BH}.

The late 1970's was also a mathematically very exciting period in other respects.
It was at this time that Atiyah and others developed new links between mathematical physics, notably quantum field
theory, on the one side, and algebraic and differential geometry on the other side. The first prime example
of this is  {\em non-abelian gauge theory} and its connections with vector bundles.
Yang-Mills fields can be described as anti-selfdual connections on $SU(2)$-bundles. Imposing asymptotic conditions on these connections at infinity, one obtains self-dual connections
on $SU(2)$-bundles on the $4$-sphere $S^4$. At this point Penrose's {\em twistor theory} comes into play. The twistor space of the $4$-sphere $S^4$ is the complex projective
$3$-space and the twistor fibration becomes $p: \mathbb P^3 \to S^4$.  Another interpretation of this map is that $\mathbb P^3=  \mathbb P(\mathbb C^4)$ is the
complex projective space and that $S^4= \mathbb P(\mathbb H^2)$ is the quaternionic projective line. Identifying $\mathbb C^4\cong \mathbb H^2$, the twistor map
$$
p:  \mathbb P^3= \mathbb P(\mathbb C^4) \to S^4= \mathbb P(\mathbb H^2)
$$
is then given by associating to a complex line in $\mathbb C^4$ the quaternionic line containing it.
Left-multiplication by $j$ moreover gives $\mathbb P^3$ a real structure $\sigma: \mathbb P^3 \to \mathbb P^3$ with $\sigma^2= - {\mathrm {id}}$.  The fibres of the twistor map
$p: \mathbb P^3 \to S^4$ are projective lines which are  invariant under $\sigma$, these are the so-called {\em {real lines}}.
Now, given a pair $(F,d)$, consisting of an $SU(2)$-bundle $F$,  together with a connection $d$, one can, by the Atiyah-Ward correspondence \cite{AW},
consider the pullback $p^*(F)$. Its associated vector bundle $E$  is a $\mathbb C^2$-bundle. The connection $d$
defines  an almost complex structure on $E$ which, since $d$ is assumed to be anti-selfdual, turns out to be integrable. In other words, $E$ is a holomorphic
 (and by GAGA thus algebraic)
rank $2$ vector bundle
on $\mathbb P^3$. Moreover, $E$ carries a real structure and, in particular, $E$ is trivial on all real lines. It is then not hard to see that $E$ is stable.
Thus one can identify self-dual $SU(2)$-connections on
$S^4$ with certain stable algebraic rank $2$  vector bundles on $\mathbb P^3$ with a real structure, the so-called {\em instanton bundles}.
Under this construction, the instanton number of the connection becomes the second Chern class of $E$. These vector bundles satisfy the additional property that the cohomology
groups $H^1(\mathbb P^3, E(-2))$ vanish. This is a translation of the fact that certain differential equations admit only trivial solutions. Now every stable rank $2$-bundle $E$ on $\mathbb P^3$ with this additional property (these are the so-called {\em mathematical instanton bundles}) are
given by  a result of Barth and Hulek \cite[Section 7]{BH} by a monad of the form
$$
n{\mathcal O}_{\mathbb P^2}(-1) \stackrel{a}  \longrightarrow (2n+2) {\mathcal O}_{\mathbb P^2} \stackrel{{}^ta\circ J}  \longrightarrow n{\mathcal O}_{\mathbb P^2}(1)
$$
where $J$ is the standard symplectic form on the vector space $\mathbb C^{2n+2}$.
The map $a$ is then nothing but an $n \times (2n+2)$-matrix with linear entries in the homogeneous coordinates of the projective plane. This can be viewed as a triple of matrices, and is often also
called  a {\em  Kronecker module}. The famous   Atiyah-Drinfeld-Hitchin-Manin (ADHM) correspondence \cite{ADHM} now says that all  instanton bundles are given by monads of the above form satisfying a reality condition.

We have already discussed that, conjecturally, all rank $2$ bundles on $\mathbb P^N$ split for $N \geq 6$. Whereas there are plenty of indecomposable rank $2$ vector bundles on
$\mathbb P^2$ and $\mathbb P^3$, only very few examples are known for higher-dimensional projective space. Essentially the only known example in characteristic $0$ is
the famous {\em Horrocks-Mumford bundle}\label{p:HM} $F=F_{{\mathrm {HM}}}$ on $\mathbb P^4$.  This is a stable rank $2$ vector bundle on $\mathbb P^4$ with Chern classes
$c_1(F)=5$ and $c_2(F)=10$. The Horrocks-Mumford bundle is a beautiful mathematical object which is distinguished by its symmetry group $N_5$ of order $15.000$. The group $N_5$ is the
semi-direct product of the Heisenberg group of level $5$, a group of order $125$, and the binary icosahedral group $SL(2,\mathbb Z/5\mathbb  Z)$. This bundle was constructed by Horrocks and Mumford in \cite{HM}
by means of a monad, but it can also be obtained via the  Serre construction from an abelian surface $A$ embedded as a surface of degree $10$ into $\mathbb P^4$.
Barth, together with Hulek and Moore, studied the many beautiful aspects of the geometry of this vector bundle in detail. The paper \cite{BHM} contains a complete classification
of all {\em Horrocks-Mumford surfaces}, i.e.~all surfaces obtained as zero-sets $A_s=\{s=0\} \subset \mathbb P^4$ of sections $0\neq s \in H^0(\mathbb P^4, F)$. These
surfaces are $(1,5)$-polarized abelian surfaces and their degenerations and this establishes a connection with (modular)  compactifications of moduli spaces of abelian surfaces.


\subsection{Algebraic surfaces}

   Barth's interest in algebraic surfaces originated from his time
   at Leiden and  became the major focus of his work since
   the 1980's. Three main lines of research can be identified in
   his work in this area:

   \begin{itemize}\compact
   \item
      fundamental work in the theory of algebraic surfaces
      \cite{AngBar82}, \cite{BarPet83}, \cite{BPV}, \cite{BHPV}, \cite{Barth-Knoerrer}
   \item
      work on abelian surfaces and Kummer surfaces  \cite{Bar85},  \cite{Bar87},  \cite{Bar93},  \cite{Barth-Bauer:conics},  \cite{Barth-Nieto:16-skew-lines},  \cite{BV} motivated originally by their
      connection with the Horrocks Mumford bundle \cite{Bar79}, \cite{BHM}

   \item
      surfaces with many symmetries and with special geometry:
      see Section~\ref{sec:symmetries}, \cite{BJAG96}, \cite{B2000}, \cite{Beven}, \cite{Barth-Bauer:conics},
      \cite{B2005}, \cite{B2007}, \cite{BS}.
   \end{itemize}

\paragraph{Algebraic surfaces -- The classification.}
   One of the central tasks of any mathematical theory is the classification of its objects.
   A systematic study of algebraic surfaces had been initialized by Max Noether, who
   as Wolf Barth, was an Ordinarius at Universit\"at Erlangen. The classification of complex
   projective surfaces was accomplished by Federigo Enriques. His works culminated
   in the book \emph{Le superficie algebriche} \cite{Enr}
   published in 1949, three years after his death.
   The classification was extended to non-algebraic compact surfaces by Kunihiko Kodaira
   in the late 60's. One of the key points in the approach of Kodaira is the study
   of elliptic fibrations of compact surfaces. Barth, together with Gerhard Angerm\"uller,
   presented in \cite{AngBar82} a complete classification of singular elliptic fibres
   on Enriques surfaces. As is typical for Barth's research, this paper contains a considerable
   list of explicit examples of Enriques surfaces and elliptic fibrations, in which fibres
   from the classification appear.

   In the paper \cite{BarPet83} with Chris Peters, Barth continued his study of
   Enriques surfaces turning attention to their automorphisms. The main result
   of this article was quite surprising: the authors proved that the automorphism
   group of a generic (in the sense of moduli) Enriques surface is large, in particular infinite,
   whereas it can be small, in particular finite, for special Enriques surfaces.
   This is quite non-intuitive, as usually, for example for curves or surfaces of general type,
   the picture is opposite: Large automorphism groups are attached to special, hence rare,
   varieties. The article, based on the global Torelli theorem for projective
   K3 surfaces, implies, as a byproduct, a result which Barth surely considered \emph{amusing}.
   It says that a generic Enriques surface has
   \begin{itemize}\compact
   \item $527$ elliptic fibrations
   \item $67\,456$ realisations as a double plane in $\P^4$
   \item $5\,396\,480$ realisations as a sextic surface in $\P^3$ passing doubly through
   the edges of a tetrahedron and
   \item $25\,903\,104$ ways to be written  as a surface of degree $10$ in $\P^5$.
   \end{itemize}
   Of course, such sample results, amusing as they might be, are only the tip of the iceberg of
   Barth's main contribution in the 80's, namely his book \emph{Compact Complex Surfaces}
   \cite{BPV} co-authored with Chris Peters and Antonius Van de Ven. This book is
   the first instance where the aforementioned Enriques-Kodaira classification
   appears in full details in print in one place. The book contains numerous, up to date at that time,
   results on surfaces of general type and on K3 and Enriques surfaces, including the global
   Torelli theorem and the theory of periods. The authors' approach to the classification theory,
   based on  Iitaka's $C_{2,1}$-conjecture, incidentally reproved in the book,
   was new and original. The book went to print in 1984. Around the same time Shigefumi Mori
   introduced a completely new way of classifying algebraic varieties of higher dimensions, known
   nowadays as the Minimal Model Programme.
   Although aimed at higher-dimensional  birational geometry, the minimal program also opens up a new view on surface
   classification. This is reflected in
   the second edition of the book, co-authored by Hulek, which appeared
   in 2004 \cite{BHPV}. It is considerably enlarged and reflects the developments of two decades, further
   including Reider's results and the ensuing improvement of the treatment of pluricanonical maps, as well as the theory of Donaldson and Seiberg-Witten invariants. The book serves still
   as a standard text for compact complex surfaces.

\paragraph{At the roots of algebraic geometry -- Equations defining algebraic varieties.}
   Barth's study of abelian varieties was strongly influenced by the series of articles
   by David Mumford
   \cite{Mum1,Mum2,Mum3} on equations defining abelian varieties. Mumford's approach
   is highly abstract. In fact, his first article contains only one equation (expressed in homogeneous coordinates):
   that of an elliptic curve
   in a Hesse pencil. There is no doubt that Barth easily handled and created complex abstract objects and arguments
   in algebraic geometry. In particular in the later stages of his career, however, his research was guided by the desire to be as explicit as possible.
   He was intrigued by the question how abstract algebraic varieties can be explicitly described
   in terms of equations involving homogeneous coordinates. The fascination by specific equations was growing
   over the years, culminating in the beautiful symmetric constructions described in Section \ref{sec:symmetries}.

   Working still on the Horrocks-Mumford bundle, Barth suggested in \cite{Bar79} a way to obtain new rank $2$
   stable algebraic vector bundles on $\P^4$. He showed that for a generic point $P\in\P^4$
   the jumping lines (see Section \ref{subsec:vectorbundles}) of $\HM$ passing through $P$
   generate a cone with vertex at $P$ over a smooth curve, which is the contact curve
   of two Kummer surfaces. Barth suggested to reverse the process and to construct vector bundles
   starting with suitable contact curves. The idea was pursued, among others, by
   Decker, Narasimhan and Schreyer, but despite these efforts it remains in the legacy
   of the not yet completed projects envisioned by Wolf Barth.

   Adler and van Moerbeke studied in \cite{AdlMoe82} algebraically integrable geodesic
   flows on $SO(4)$ and related them to affine parts of abelian surfaces in $\C^6$
   defined as complete intersections of $4$ quadrics. These affine parts come from
   abelian surfaces embedded in $\P^7$ by complete linear systems of type $(2,4)$.
   Inspired by these results, Barth studied in \cite{Bar85} abelian surfaces
   with $(1,2)$ polarization and gave a complete description of these surfaces and  a description of their moduli space.
   In a paper dedicated to Friedrich Hirzebruch \cite{Bar87}, Barth studied the question
   whether there exist other abelian surfaces, besides those discovered by Adler and van Moerbeke,
   which might be related to integrable Hamiltonian systems. He studied explicit
   equations of such potential abelian surfaces and their symmetries imposed by Heisenberg
   groups and arrived at the conclusion that
   among abelian surfaces with general moduli
   no such surfaces exist.

   Wolf Barth continued his studies of abelian surfaces defined by quadratic equations in
   \cite{Bar93}. There he considered principally polarized abelian surfaces $A$
   embedded in $\P^8$ by the third power of the theta divisor. By a result of Kempf \cite{Kempf} it was known that the homogeneous ideal
   of such surfaces is generated by forms of degree $2$ and $3$.
   Barth showed that quadrics suffice to generate the {\em ideal sheaf} if and only if the
   polarized abelian surface $A$ is not a product of two elliptic curves.
   This paper is yet another example of Barth's interest in concrete
   equations. He provides explicit equations of quadrics cutting out $A\subset\P^8$.

   If $A$ is an abelian surface, then its quotient $X=A/{(-1)_A}$ by the $(-1)$-involution
   is by definition a \emph{Kummer surface}. Its desingularization $\tilde X$, the \emph{smooth Kummer surface}
   of $A$, is then a $K3$ surface which contains 16 skew smooth rational curves
   (corresponding to the half-periods of $A$).
   It is a result of Nikulin \cite{Nikulin:Kummer}
   that, conversely, every K\"ahler K3 surface with 16 skew smooth rational curves
   arises in this way as a smooth Kummer surface.
   Barth
   was interested in finding such surfaces with many skew rational curves
   in \emph{three-space} -- always with an eye towards concrete realizations of surfaces.
   In joint work with Nieto he
   rediscovered in \cite{Barth-Nieto:16-skew-lines}
   the smooth Kummer surfaces found by Traynard \cite{Traynard}. These contain 16 skew lines
   and they are associated with abelian surfaces of type $(1,3)$.
   It was his idea that one should be able to generalize the construction in order
   to obtain also smooth Kummer surfaces in $\P^3$ with 16 skew \emph{conics}
   (instead of lines).
   This was in fact possible and
   the surfaces constructed in this way turn out to be very interesting
   as they
   contain a surprisingly high number
   of conics in total \cite{Barth-Bauer:conics}.
   It was shown later \cite{Bauer:smooth-Kummer} that it is
   even possible to obtain
   16 skew smooth rational curves of \emph{any} given degree.

\paragraph{Visualizations of Algebraic surfaces -- From
   the first steps to the Imaginary Exhibition.}

   In addition to his theoretical work, Barth had a great
   interest in visualizations of algebraic surfaces.
   He had a fascination for the
   models of algebraic surfaces that were produced in the
   19th century (made in plaster, sometimes in
   wood), and he
   loved them both
   for the geometric insight they could provide
   as well as from a purely esthetical point of view.
   His contribution (with Horst Knörrer) to the volume
   \emph{Mathematical Models}~\cite{Barth-Knoerrer}
   shows his
   ample knowledge and his
   appreciation of these kinds of visualizations.

   With the availability of ever more powerful computers,
   it was his idea that
   given the polynomial equation $f(x,y,z)=0$ of an algebraic surface
   in three-space, it should be possible
   to generate convincing images of the (real) surface (i.e.~of the zero set
   of $f$ in $\R^3$).
   The crucial breakthrough in this direction was
   achieved in the diploma thesis by
   Stephan Endrass (1992), supervised by Barth.
   In this thesis,
   the program \texttt{surf} was developed (on a Cadmus
   workstation),
   which was able to produce
   impressive images of algebraic surfaces from their equations.
   When Barth had constructed his now famous sextic surface~\cite{BJAG96}
   (see Section \ref{sec:symmetries}),
   it was a great moment for him to actually \emph{see} the
   surface after he had found it by theoretical means.
   The program \texttt{surf} was extended and transferred to other systems
   during the 90s in Barth's group in Erlangen.
   Since 2000, further development has been done by the algebraic geometry group
   in Mainz,
   and the program
   is still in constant use nowadays
   via the GUI program \texttt{surfer}
   (\texttt{https://imaginary.org/de/program/surfer}),
   which
   relies on the \texttt{surf} kernel for its computations.
   So we can enjoy the
   impressive images of algebraic surfaces
   in the Imaginary project
   (\verb|https://imaginary.org|)
   thanks to
   Barth's early interest in visualizations.
   Of course, Barth's sextic surface
   deserves its prominent place in the exhibition.


\subsection{Symmetries}\label{sec:symmetries}

\begin{figure}[t]
\minipage{0.40\textwidth}
\begin{center}
\includegraphics[width=5cm, height=5cm]{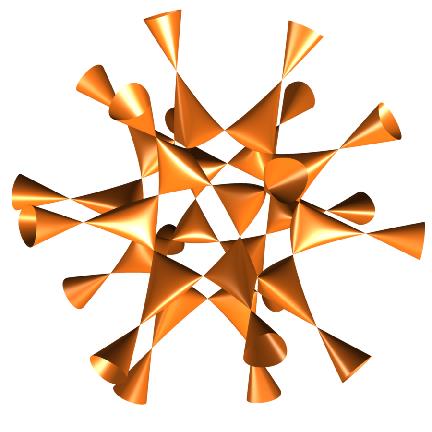}
\caption{Barth's sextic surface}
\label{barth_sextic}
\end{center}
\endminipage\hfill
\minipage{0.40\textwidth}
\begin{center}
\includegraphics[width=5cm, height=5cm]{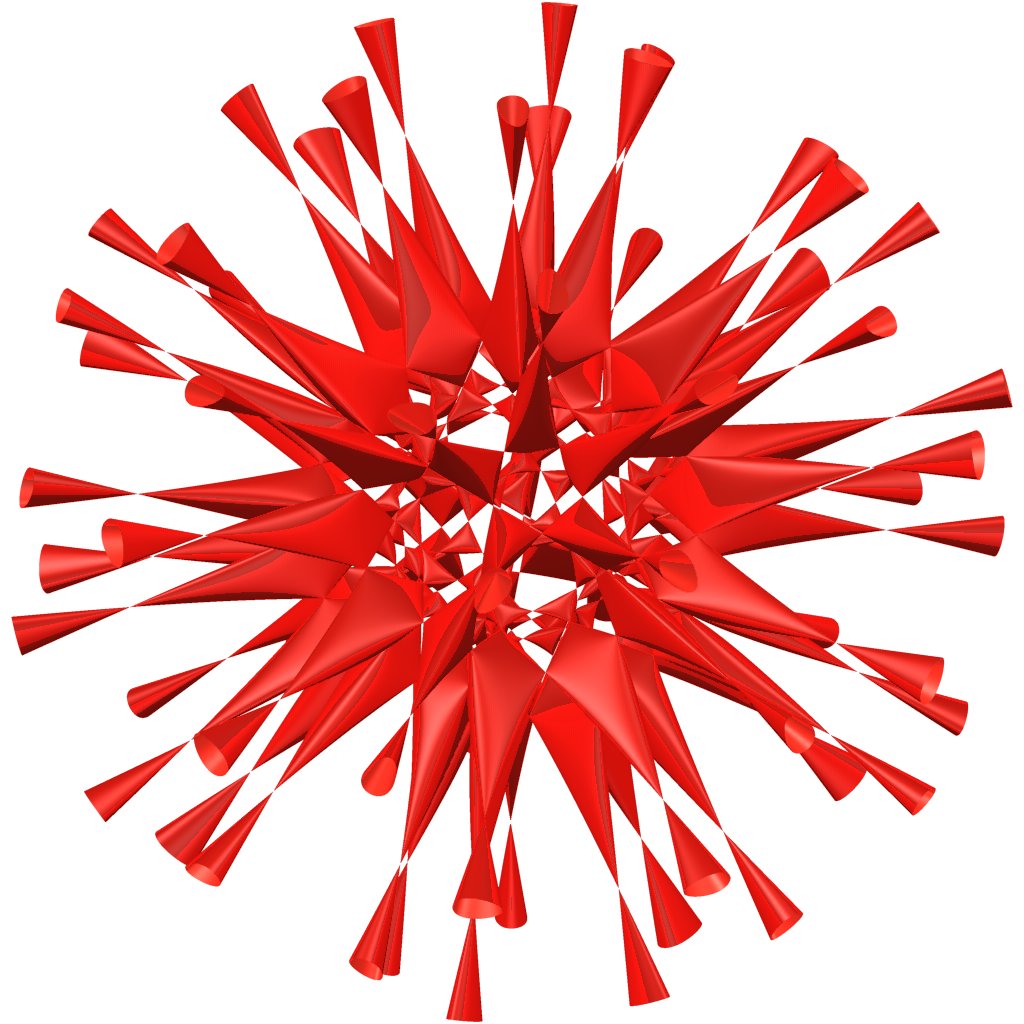}
\caption{Barth's decic surface}
\label{decic}
\end{center}
\endminipage
\end{figure}

In the last 20 years of his research activity Wolf Barth
got very much interested in the use of
symmetries in algebraic geometry. He was fascinated by the beauty of the surfaces
that one can produce using symmetries. In particular, he investigated
surfaces with
many rational curves (resp.~singularities obtained by contracting  configurations of such curves) and constructed numerous examples.

Recall that an $A_1$ (resp.~an $A_2$) point of a surface is a singularity locally given by the equation
$$
x \cdot y - z^2 = 0
$$
(resp.~$x \cdot y - z^3 = 0$). Such a point can be locally obtained as the quotient of the two-dimensional complex unit ball
by an appropriate
$\Z/2\Z$ (resp.~$\Z/3\Z$) action. An $A_1$ (resp.~$A_2$) singularity is also called a {\it node} (resp.~a {\it cusp}).

A foundational work of W.~Barth on this subject is the paper
of 1996, \cite{BJAG96},
where he discovered what is now called {\it Barth's sextic} (see Figure \ref{barth_sextic}).
It is a surface of degree 6 in three dimensional complex
projective space with 65 nodes.

In this paper he considers the symmetry group $I$ of the icosahedron  in euclidean three dimensional space $\R^3$.
It is well-known that this is isomorphic to the alternating group $A_5$ and it acts on the ring of coordinates $\R[x,y,z]$.
The ring of invariant polynomials
is generated by polynomials of degree 2, 6 and 10 respectively. The polynomial of degree 2  can be written as $x^2+y^2+z^2$ and the other two invariants
of degree 6 and 10  can be written in the following forms (they were already known to Goursat \cite{Go}):
\be
   & Q(x,y,z)=(\tau^2x^2-y^2)(\tau^2y^2-z^2)(\tau^2z^2-x^2) \\
   & R(x,y,z)=(x^2-\tau^4y^2)(y^2-\tau^4z^2)(z^2-\tau^4x^2)(x+y+z)(x+y-z)(x-y+z)(x-y-z)
\ee
where $\tau=(1+\sqrt{5})/2$ is the golden ratio. Combining these two equations with the equation of the three dimensional sphere $S: x^2+y^2+z^2-1=0$, Barth obtains two families of
symmetric surfaces invariant under the action of $I$. Barth's sextic  belongs to the family of surfaces
with equation in affine coordinates
$$
f_{\alpha}:\,Q(x,y,z)-\alpha (x^2+y^2+z^2-1)^2=0.
$$
For generic choice of the complex parameter $\alpha$ the surface has 45 nodes that belong to special lines
of the icosahedron, the {\it mid lines}, as Barth called them, that connect two opposite mid points of edges. Then Barth imposes geometric
conditions in a clever way to get an extra orbit of 20 nodes for the action of $I$. He finally  obtains
that by taking the parameter $\alpha=(2\tau+1)/4$ the
surface has in total 65 nodes. A very similar construction was then used by Barth to combine the degree ten polynomial $R(x,y,z)$ and the
equation of the sphere $S$ to obtain {\it Barth's decic} (Figure \ref{decic}) that has 345 nodes. As in the previous case, the singular points are located on special
families of lines and planes related to the geometric properties of the icosahedron (the 10 planes through the origin, parallel to the twenty faces of the icosahedron).

In the paper \cite{M} of 1984, Miyaoka gave a bound for the maximum number of nodes that a surface of given degree in complex projective three-space can have.
The bound is $66$ for a surface of degree $6$ and $360$ for a surface of degree $10$. The two examples produced by Barth
are the best known examples so far in these two degrees. In particular
 the discovery of Barth's sextic was quite surprising, since it was thought that
no surface of degree $6$ surface with $65$ nodes could exist:
Fifteen years before, Catanese and Ceresa  \cite{CC} had erroneously claimed
that $64$ was the maximum number of ordinary double points for a sextic surface. Soon after Barth's discovery, Jaffe and Ruberman \cite{JR} showed, using coding theory, that in fact $65$ is the maximum possible number of nodes for surfaces of degree $6$.
The problem of the maximum number of nodes in degree six was thus solved!

Barth's sextic  remains until now an exceptional and beautiful example
of how symmetries are a powerful tool in attacking problems in algebraic geometry.

The lower bound of 345  given by Barth's decic  also remains to our knowledge the best so far. It is still not known whether Miyaoka's bound of 360 nodes is possible or not for a surface of degree 10 in complex three-dimensional projective space.
The idea of Barth of considering surfaces with many symmetries to attack the problem of finding surfaces with many nodes, paved the way to the discovery of more {\it world record surfaces}, as Barth termed them, and  	greatly influenced several mathematicians working on the subject, in particular Barth's PhD students: S. Endrass and A. Sarti. They found respectively a surface of degree 8 with 168 nodes in 1997,  \cite{End}, and a surface of degree 12 with 600 nodes in 2001, \cite{S}. These numbers of nodes approach the bounds of Miyaoka which are respectively 174 and 645 and are so far the best known examples. Both the surfaces are very symmetric, they have respectively the symmetries of some extension of order two of the dihedral group $D_8$ and of the bipolyhedral icosahedral group of order 7200.

The ideas of Barth affected the works of several other mathematicians such as D.~van Straten and his students in Mainz, as well as S. Cynk in Crakow. His ideas also influenced the study of other difficult and classical problems in algebraic geometry such as the study of the maximum number of lines on projective surfaces. Around 2001 Barth found a smooth quintic surface with 75 lines that was later described in the paper \cite{Xie}
(the author, as he said, followed a suggestion by W.~Barth).

In a paper of 2003 W.~Barth and A.~Sarti \cite{BS} studied the relation of the symmetric
surfaces of the PhD thesis of A. Sarti, which had been conducted under the supervision of W.~Barth, and
 K3 surfaces. In particular, they show that the minimal resolution of the quotient of the surface with 600 nodes by the bipolyhedral icosahedral group is a K3 surface with maximal Picard number, namely 20. Wolf Barth had always had a special interest in the beauty of the geometry of K3 surfaces and in his paper \cite{Beven} he studied divisible sets of rational curves on K3 surfaces. More precisely, let $L_1, \ldots, L_k$, $k\geq 1$, be smooth disjoint rational curves on a K3 surface,
Nikulin in \cite{Nik}
 showed that if this set is {\it 2-divisible}, i.e.~the sum $L_1+\ldots+L_k$ is equivalent to two times a divisor in the N\`eron-Severi group,
 then $k=8$ or $k=16$ and $16$ disjoint rational curves are always an even set. This is not the case for eight disjoint rational curves. In
\cite{Beven} Barth  characterizes even sets of 8 rational curves on the most geometric projective models of K3 surfaces, such as double covers of $\IP^2$ ramified along a sextic curve, quartics in $\IP^3$ and double covers of $\IP^1\times\IP^1$ ramified along a curve of bi-degree $(4,4)$. The results provide a geometric expression of the lattice-theoretical results of Nikulin and they deeply affected the work of several other mathematicians, such as B.~van Geemen and his students in Milan, working in particular on symplectic involutions on K3 surfaces (which are strictly related to the existence of such even sets).

Wolf Barth also investigated  sets of $A_2$ singularities on surfaces.
In general,  techniques of enumerative geometry are not sensitive enough to give precise statements on the number of A-D-E points on degree $d$
surfaces in $\IP^3$. In the case of quintics this problem was circumvented by Beauville, who associated to a surface with  $\mu$ nodes  a linear code in $\F_2^{\mu}$.
Each word of the code in question is given by a so-called {\em even set of nodes}.  By definition, a set of nodes $P_1, \ldots, P_k$ on a surface
$X$ is even if and only if the local $\Z/2\Z$-quotient structure around the $k$ nodes can be defined globally, i.e.~when there exists a double cover
$Y$ of $X$ branched only at the points $P_1, \ldots, P_k$ with $Y$ smooth along the  branch locus.
Once the code is defined, one can use  coding theory to
arrive at bounds on the number of singularities and constraints on their configuration.
In the papers \cite{B1998}, \cite{B2000} Barth  generalized the above notion to surfaces with $A_2$ singularities.
In this case one arrives at a ternary code. Barth used this generalization
to show that
a cuspidal quartic surface can carry at most eight singularities and
every K3-surface with nine $A_2$ points arises as a 3:1 quotient of a complex torus.
This result generalizes the classical and well known construction of Kummer surfaces.
Furthermore, he classified all tori occurring in the construction
and gave explicit examples of quartics with eight $A_2$ points. Finally, the papers \cite{B2005}, \cite{B2007}, joint with S.~Rams, contain sharp effective bounds on
the minimal weights of ternary codes given by low-degree surfaces and the computation of codes for certain surfaces.

It should be pointed out that Barth's work on divisible sets of singularities played an important role in the  project of
classifying fundamental groups of open Enriques and K3 surfaces which has recently been carried out by J.H.~Keum and D.Q.~Zhang,
whereas his published and unpublished examples appear in various contexts in algebraic geometry (e.g. Barth's quintics surface with 75 lines turns out to be the smooth quintic with the highest Picard number known so far -- see \cite{RS41}).


\section{Wolf Barth -- The teacher}

   Over the years, Wolf Barth developed an ever growing interest
   in teaching mathematics, which became apparent in a number of
   activities that are non-routine among top research
   mathematicians: For instance, when he found that the mathematics education
   regularly provided for elementary and middle school teachers
   was not up to the standards that he aimed at, rather than
   theorizing about this fact, or putting blame on others, he
   volunteered to restructure the courses in question and to
   teach them himself for a number of years to come -- an
   instance of the
   hands-on approach he employed in such situations.

   As far as standard courses such as Analysis, Linear Algebra, Abstract
   Algebra and Complex Analysis are concerned, he developed
   his unique way of attacking the subjects,
   enriched by
   written
   manuscripts for the students which were not just compilations
   of the known textbook approaches. One such example is his
   \emph{Linear Algebra} text, which later became a Springer book
   coauthored by Peter Knaber \cite{Barth-Knabner}.

   Also, it was his intention to convey to students his
   fascination for classical topics in geometry. As an example,
   he developed and taught a course on \emph{circles}
   \cite{Barth:circles}.
   In his motivation one can see parallels with his activities
   in the environment protection movement: Barth knew that classical
   geometry lost its place at universities because it was considered
   by the majority of academia mathematicians as too elementary.
   Barth himself considered numerous beautiful theorems related to circles as members of a rare
   species which deserved protection from vanishing from mankind's intellectual
   legacy. His attitude to projective geometry was very similar.

   Over the years Barth had been continuously involved in the training of
   the next generation of algebraic geometers. His first student, Wilfred Hulsbergen,
   obtained his PhD in 1976, the last student of Barth, Thomas Werner,
   graduated in 2012. Four of the authors of this obituary were Barth's students
   and the remaining one, Rams, obtained his habilitation with Barth in Erlangen.

   The last of the named authors came from Cracow to Erlangen, exactly on the day of
   German reunification, as an exchange student. He encountered Barth and became
   his graduate student. It was the time of rapid political changes and growing
   hopes throughout Europe. All of a sudden, Eastern Europe became part of the Free World.
   Direct scientific exchange of ideas and people was possible and it was generously
   supported by various organizations including DFG.
   Barth successfully applied for considerable grants supporting the library
   of the Institute of Mathematics of the Jagiellonian University and for four years
   he coordinated a staff exchange program with Cracow. He visited Cracow himself twice.
   His visits and the visits of Polish young mathematicians, including S\l awomir Cynk,
   Zbigniew Jelonek and Piotr Tworzewski resulted in transplanting modern algebraic
   geometry to Cracow, which nowadays is, next to Warsaw, the strongest centre
   of algebraic geometry in Poland. This development would have been hardly possible without
   Barth's support and engagement.

Barth also had important connections to Italy.
At the end of the 1990's, Graziano Gentili,
who at that time was the Head of the {\it Scuola Matematica Interuniversitaria}
invited Wolf Barth to teach a course in algebraic geometry at the Perugia summer school.
This is a very famous Italian summer school aiming to prepare young Italian and foreign
students to work on a PhD thesis.
Barth gave the algebraic geometry course both  in 1997
and in 2000, and his lectures were highly appreciated.
The fourth named author
met Barth at the school in 1997 and then started a PhD thesis in Erlangen with him.
Over the years Barth became a good friend of Gentili, whom he  invited several times to Erlangen.

\begin{center}
   * * *
\end{center}

Wolf Barth had a lively personality. He was very much interested
in politics, history and contemporary events. He rejected any form of xenophobia
and he was actively involved in the protection of the environment. He
loved every kind of beauty (flowers, art, music, and all geometric
shapes). Material
things were not important to him. He loved his family, his children and he was
loved in return.
He was a brilliant mathematician, passionate for algebraic geometry
where he obtained  fundamental results and wrote a foundational book.
His famous sextic surface, the Barth surface, is one of the icons of algebraic surfaces.
To all who came in contact with him he communicated his great passion and enthusiasm for mathematics.

\paragraph{Acknowledgements.}
   We warmly thank Matthias Barth, Renate
   Ararat, Herbert Lange, Laurent Gruson, Chris Peters, Ulf Persson,
   and the administration of
    Universit\"at Erlangen
   for their valuable help during the preparation of the manuscript.



\begin{center}
   \includegraphics[width=\linewidth]{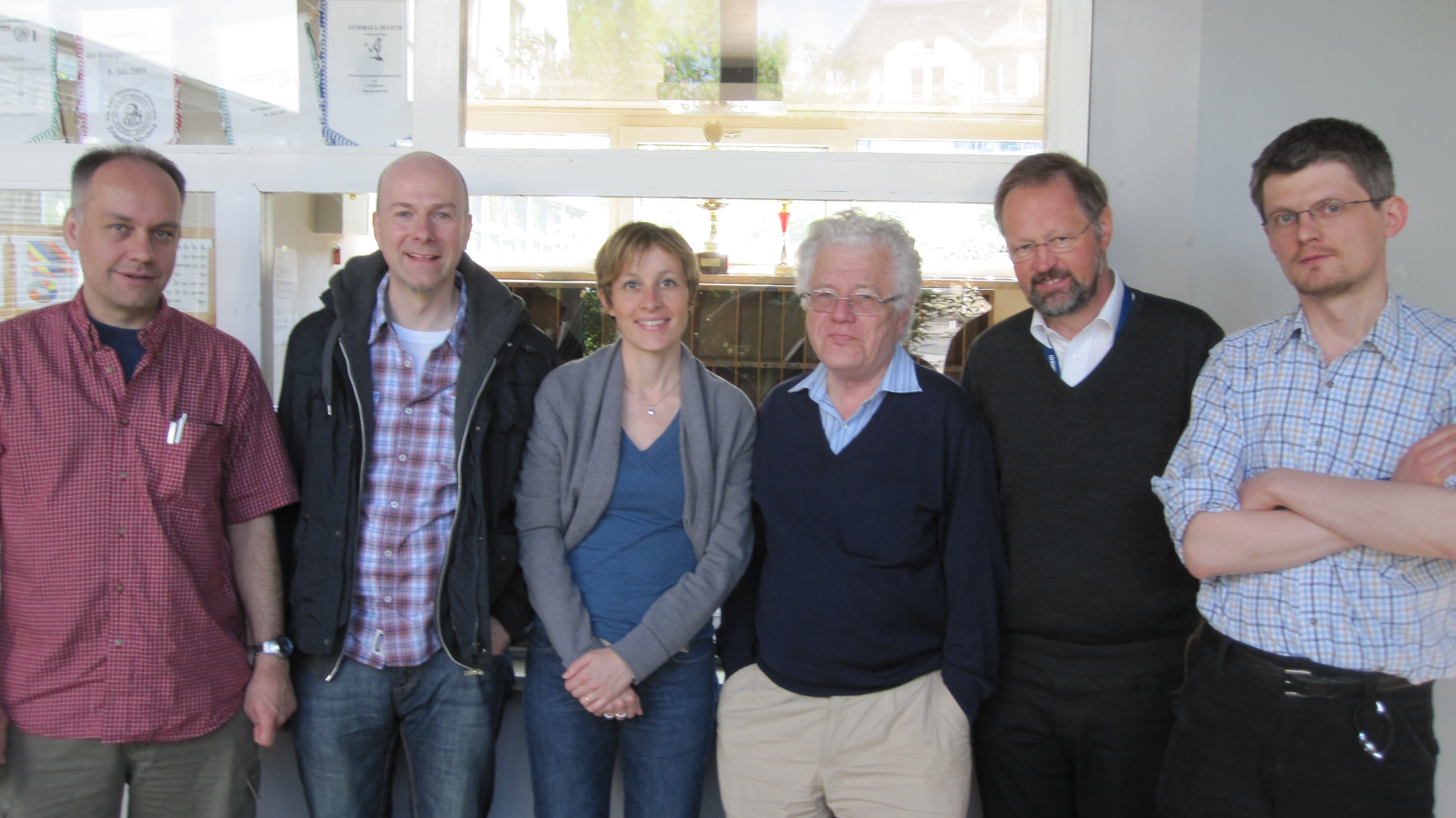}
   Prof.~Dr.~Wolf Barth and his students
   at the conference \emph{Groups and Algebraic Geometry},
   Erlangen, April 2011.
   (From left to right: T.~Szemberg, Th.~Bauer, A.~Sarti, W.~Barth, K.~Hulek, S.~Rams)
\end{center}



\begin{thebibliography}{99}\footnotesize\compact\frenchspacing

\bibitem{AdlMoe82}
   Adler, M., and van Moerbeke, P.:
   The algebraic integrability of geodesic flow on SO(4).
   Invent. Math. \textbf{67}, 297-331 (1982)

\bibitem{AngBar82}
   Angerm\"uller, G., and Barth, W.:
   Elliptic fibres on Enriques surfaces.
   Compositio Math. {\bf 47}, 317-332 (1982).

\bibitem{ADHM}
   Atiyah, M.~F., Hitchin, N.~J., Drinfeld, V.~G., and Manin, Yu.~I.:
   Construction of instantons.
   Phys. Lett. A. {\bf 65}, 185--187 (1978).

 \bibitem{AW}
   Atiyah, M.~F., and Ward. R.~S.:
   Instantons and algebraic geometry.
   Comm. Math. Phys. {\bf 55}, 117--124 (1977).

\bibitem{B1}
   Barth, W.:
   Fortsetzung, meromorpher {F}unktionen in {T}ori und komplex-projektiven {R}\"aumen.
   Invent. Math. {\bf 5}, 42--62 (1968).

\bibitem{B2}
   Barth, W.:
   Transplanting cohomology classes in complex-projective space.
   Amer. J. Math. {\bf 92}, 951--967 (1970).

\bibitem{B3}
    Barth, W.:
    Some properties of stable rank-2 vector bundles on $\mathbb P_n$.
    Math. Ann. {\bf 226}, 125--150 (1977).

 \bibitem{B4}
    Barth, W.:
    Moduli of vector bundles on the projective plane.
    Invent. Math. {\bf 42}, 63--91 (1977).

\bibitem{Bar79}
   Barth, W.:
   Kummer surfaces associated with the Horrocks-Mumford bundle.
   Journ\'ees de G\'eometrie Alg\'ebrique d'Angers, Juillet 1979/Algebraic Geometry, Angers, 1979, pp. 29--48,
   Sijthoff \& Noordhoff, Alphen aan den Rijn--Germantown, Md., 1980.

\bibitem{Bar85}
   Barth, W.:
   Abelian surfaces with $(1,2)$--polarization.
   Algebraic geometry, Sendai, 1985, 41--84, Adv. Stud. Pure Math., 10, North-Holland, Amsterdam, 1987.

\bibitem{Bar87}
   Barth, W.:
   Affine parts of abelian surfaces as complete intersections of four quadrics.
   Math. Ann. \textbf{278}, 117-131 (1987).

\bibitem{Bar93}
   Barth, W.:
   Quadratic equations for level--3 abelian surfaces.
   Abelian varieties (Egloffstein, 1993), 1--18, de Gruyter, Berlin, 1995.

\bibitem{BJAG96}
   Barth, W. :
   Two projective surfaces with many nodes, admitting the symmetries of the icosahedron.
   J. Algebraic Geom., {\bf 5}, 173--186 (1996).

\bibitem{B1998}
    Barth, W. :
   K3 surfaces with Nine Cusps. Geometriae Dedicata {\bf 72}: 171--178, (1998).

\bibitem{B2000}
   Barth, W. :
   On the Classification of K3 Surfaces with Nine Cusps.
   Complex Analysis and Algebraic Geometry -- A Volume in Memory of Michael Schneider. Editors: T. Peternell, F.-O. Schreyer, 42--59 (2000).

\bibitem{Beven}
   Barth, W.:
   Even sets of eight rational curves on a {$K3$}-surface,
   Complex Geometry ({G}\"ottingen, 2000), 1--25, Springer, Berlin, (2002).

\bibitem{Barth:circles}
   Barth, W.:
   Kreise. Skriptum, Erlangen 1997.

   \verb|https://www.studium.math.fau.de/fileadmin/studium/skripten/barth/kreiset1.ps|

\bibitem{BL}
   Barth, W., and Larsen, M.~E.:
   On the homotopy groups of complex projective algebraic  manifolds.
   Math. Scand. {\bf 30}, 88--94 (1972).

\bibitem{BH}
   Barth, W., and Hulek, K.:
   Monads and moduli of vector bundles.
   Manuscripta Math. {\bf 25}, 323--347 (1978).

\bibitem{BHM}
   Barth, W., Hulek, K., and Moore, R.:
   Degenerations of {H}orrocks-{M}umford surfaces.
   Math. Ann. {\bf 277}, 735--755 (1987).

\bibitem{BarPet83}
   Barth, W.~P., and Peters, C.~A.~M.:
   Automorphisms of Enriques surfaces.
   Invent. Math. {\bf 73}, 383--411 (1983).

\bibitem{BPV}
   Barth, W.~P., Peters, C.~A.~M., and Van de Ven, A.:
   Compact complex surfaces.
   Ergebnisse der Mathematik und ihrer Grenzgebiete. 3. Folge. Vol. 4,  Springer-Verlag, Berlin, 1984.

\bibitem{BHPV}
   Barth, W.~P., Hulek, K., Peters, C.~A.~M., and Van de Ven, A.:
   Compact complex surfaces. Second edition.
   Ergebnisse der Mathematik und ihrer Grenzgebiete. 3. Folge. Vol. 4, Springer-Verlag, Berlin, 2004.

\bibitem{Barth-Knoerrer}
   Modelle algebraischer Flächen. Kommentare von W.~Barth und H. Knörrer. Mathematische
   Modelle, G. Fischer Ed., 7-24, Vieweg (1986).

\bibitem{Barth-Nieto:16-skew-lines}
   Barth, W., and  Nieto, I.:
   Abelian surfaces of type (1,3) and quartic surfaces with 16 skew lines. J. of Alg. Geometry, Vol. 3,
   No. 2, 173-222 (1994).

\bibitem{B2005}
    Barth, W., and Rams, S.  :
    Equations of low-degree projective surfaces with three-divisible sets of cusps. Math. Z. {\bf 249},
    283--295 (2005).

   \bibitem{B2007}
    Barth, W., and Rams, S.  :
Cusps and Codes. Math. Nachr. {\bf 280}, 50--59 (2007).

\bibitem{BS}
   Barth, W., and  Sarti, A.:
     Polyhedral groups and pencils of {$K3$}-surfaces with maximal
              {P}icard number,
   Asian J. Math., {\bf 7}, 519--538, (2003).

\bibitem{BVdV}
  Barth, W., and Van de Ven, A.:
  A decomposability criterion for algebraic {$2$}-bundles on projective spaces.
  Invent. Math. {\bf 25}, 91--106 (1974).

\bibitem{BV}
 Barth, W., and Verra, A.:
   Torsion on {$K3$}-sections, Problems in the theory of surfaces and their classification.
              ({C}ortona, 1988), Sympos. Math., XXXII, 1--24, Academic Press, London, (1991).

\bibitem{Barth-Bauer:conics}
   Bauer, Th., and  Barth, W.:
   Smooth quartic surfaces with 352 conics.
   Manuscripta math. 85, 409-417 (1994)

\bibitem{Bauer:smooth-Kummer}
   Bauer, Th.:
   Smooth Kummer surfaces in projective three-space.
   {Proc. Amer. Math. Soc.} {\bf 125}, 2537-2541 (1997)

\bibitem{Be}
   Be{\u\i}linson, A. A.:
   The derived category of coherent sheaves on {${\bf P}^n$}.
   Selecta Math. Soviet. {\bf 3}, 233--237 (1983/84).

   \bibitem{CC}
    Catanese, F., and Ceresa, G. :
     Constructing sextic surfaces with a given number {$d$}\ of
              nodes.
 J. Pure Appl. Algebra, {\bf 23}, 1--12, (1982).

 \bibitem{DLP}
  Drezet, J.-M., and Le Potier, J.:
  Conditions d'existence des fibres stables de rang \'elev\'e sur {${\bf P}_2$}.
  Vector bundles on algebraic varieties ({B}ombay, 1984), 133--158 (1987).
  Tata Inst. Fund. Res. Stud. Math. {\bf 11}.

  \bibitem{Ell}
  Ellingsrud, G.:
  Sur l'irr\'eductibilit\'e du module des fibr\'es stables sur {${\bf P}^{2}$}
  Math. Z. {\bf 182}, 189--192 (1983).


 \bibitem{End}
   Endrass, S. :
A projective surface of degree eight with {$168$} nodes.
J. Algebraic Geom.,{\bf6}, 325--334, (1997).


\bibitem{Enr}
Enriques, F. :
     Le {S}uperficie {A}lgebriche.
Nicola Zanichelli, Bologna, 1949.

\bibitem{FHS}
  Forster, O., Hirschowitz, A., and Schneider, M.
  Type de scindage g\'en\'eralis\'e pour les fibr\'es stables.
   Vector bundles and differential equations ({P}roc. {C}onf.,  {N}ice, 1979),
   Progr. Math. {\bf 7}, 65--81 (1980).

\bibitem{FH}
   Fulton, W., and Hansen, J.:
   A connectedness theorem for projective varieties, with applications to intersections and singularities of mappings.
   Ann. of Math. (2) {\bf 1}, 159--166 (1979).

\bibitem{Gi}
   Gieseker, D.:
   On the moduli of vector bundles on an algebraic surface.
   Ann. of Math. (2) {\bf 1}, 45--60 (1977).

   \bibitem{Go}
   Goursat, E. :
   Etude des surfaces qui admettent tous les plans de symetrie d'un polyedre regulier.
   Ann. Sci. Ecole Norm. Sup. (3) {\bf 4}, 159--200 (1887).

\bibitem{GM}
   Grauert, H., and M{\"u}lich, G.:
   Vektorb\"undel vom {R}ang {$2$} \"uber dem {$n$}-dimensionalen komplex-projektiven {R}aum.
   Manuscripta Math. {\bf 16}, 75--100 (1975).


\bibitem{H}
    Hartshorne, R.:
    Varieties of small codimension in projective space.
    Bull. Amer. Math. Soc. {\bf{80}}, 1017--1032 (1974).

\bibitem{Hi}
    Hirschowitz, A.:
    Sur la restriction des faisceaux semi-stables.
    Ann. Sci. \'Ecole Norm. Sup. (4), 199--207 (1981).

\bibitem{Ho}
    Horrocks, G.:
    Vector bundles on the punctured spectrum of a local ring.
    Proc. London Math. Soc. (3) {\bf 14}, 689--713 (1964).

\bibitem{HM}
    Horrocks, G., and Mumford, D.:
    A rank {$2$} vector bundle on {${\bf P}^{4}$} with {$15,000$}\ symmetries.
    Toplology {\bf 12}, 63--81 (1973).

 \bibitem{Hul}
      Hulek, K.:
      Stable rank-{$2$} vector bundles on {${\bf P}_{2}$} with {$c_{1}$}\ odd.
      Math. Ann. {\bf 242}, 241--266 (1979).

    \bibitem{JR}
    Jaffe, D. B., and Ruberman, D. :
   A sextic surface cannot have {$66$} nodes.
  J. Algebraic Geom., {\bf 6}, 151--168 (1997).

    \bibitem{Kempf}
   Kempf, G. R.:
   Projective coordinate rings of abelian varieties,
Algebraic analysis, geometry, and number theory ({B}altimore,
              {MD}, 1988), Johns Hopkins Univ. Press, Baltimore, MD,
               225--235, (1989).

\bibitem{Barth-Knabner}
   Knabner, P., and Barth, W.:
   Lineare Algebra. Grundlagen und Anwendungen.
   Springer, 2013.

\bibitem{L}
   Larsen, M.~E.:
   On the topology of complex projective manifolds.
   Invent. Math. {\bf 19}, 251--260 (1973).

\bibitem{LP}
   Le Potier, J.:
   Fibr\'es stables de rang {$2$} sur {${\bf P}_{2}({\bf C})$}.
   Math. Ann. {\bf  241}, 217--256 (1979).


\bibitem{M1}
   Maruyama, M.:
   Moduli of stable sheaves. {I}.
   J. Math. Kyoto Univ. {\bf 3},  557--614 (1978).

\bibitem{M2}
   Maruyama, M.:
   Moduli of stable sheaves. {II}.
   J. Math. Kyoto Univ. {\bf 1},  91--126 (1977).

\bibitem{M3}
   Maruyama, M.:
   The rationality of the moduli spaces of vector bundles of rank {$2$} on {${\bf P}^2$}.
   Algebraic geometry, {S}endai, 1985, Adv. Stud. Pure Math. {\bf 10}, 399--414 (1987).

   \bibitem{M}
   Miyaoka, Y.:
     The maximal number of quotient singularities on surfaces with
              given numerical invariants.
 Math. Ann., {\bf 268}, 159--171, (1984).

\bibitem{Mum1}
   Mumford, D.:
   On the equations defining abelian varieties. I.
   Invent. Math. \textbf{1}, 287--354, (1966).

\bibitem{Mum2}
   Mumford, D.:
   On the equations defining abelian varieties. II.
   Invent. Math. \textbf{3}, 75--135, (1967).

\bibitem{Mum3}
   Mumford, D.:
   On the equations defining abelian varieties. III.
   Invent. Math. \textbf{3}, 215--244, (1967).

\bibitem{Nik}
   Nikulin, V. V.:
   Kummer surfaces.
   Izv. Akad. Nauk SSSR Ser. Mat., {\bf 39}, 278--293, (1975).

\bibitem{Nikulin:Kummer}
   Nikulin, V.V.:
   On Kummer surfaces.
   Transl. to English, Math. USSR. -- Izv. 9, 261-275 (1975)

 \bibitem{OSS}
   Okonek, C., Schneider, M., and Spindler, H.:
   Vector bundles on complex projective spaces.
   Corrected reprint of the 1988 edition,
    With an appendix by S. I. Gelfand.
    Modern Birkh\"auser Classics (2011).


 \bibitem{RS41}
Rams, S., and Sch\"utt, M.:
The Barth quintic surface has Picard number 41.
Ann. Sc. Norm. Super. Pisa Cl. Sci. (5) {\bf XIII} (2014), 533--549.

 \bibitem{S}
    Sarti, A.:
   Pencils of symmetric surfaces in {${\mathbb P}_3$},
  J. Algebra, {\bf 246}, 429--452,(2001).

 \bibitem{Sch}
     Schwarzenberger, R. L. E.:
     Vector bundles on the projective plane
     Proc. London Math. Soc. (3), {\bf 11}, 623--640, (1961).

 \bibitem{Stro}
     Str{\o}mme, S. A.:
     Deforming vector bundles on the projective plane.
     Math. Ann. {\bf 263}, 385--397 (1983).

\bibitem{Traynard}
   Traynard, E.:
   Sur les fonctions th\^eta de deux variables et les surfaces hyperelliptiques.
   Ann. scient. \'Ec. Norm. Sup., 3. s\'er., t. 24, 77-177 (1907).

  \bibitem{Xie}
    Xie, J. :
   More quintic surfaces with 75 lines.
   Rocky Mountain J. Math., {\bf 40}, 2063--2089, (2010).

\end{thebibliography}
\end{document}